\documentclass[12pt]{article}
\usepackage{amsfonts}

\usepackage[all]{xy}

\def\squarebox#1{\hbox to #1{\hfill\vbox to #1{\vfill}}}
\newcommand{\qed}{\hspace*{\fill}
\vbox{\hrule\hbox{\vrule\squarebox{.667em}\vrule}\hrule}\smallskip}
\newtheorem{teorema}{Theorem}[section]

\newenvironment{profe}{\noindent {\bf Proof:}}{\hfill $\qed $ \newline}
\begin{document}

\title{Existence of complete Lyapunov functions \\ for semiflows on separable metric spaces}
\author{Mauro Patr\~{a}o
\thanks{Department of Mathematics. UnB, Brazil. Supported by CNPq grant no.\ 310790/09-3}
}

\maketitle

\begin{abstract}
\let\thefootnote\relax
\let\thefootnote\relax\footnotetext{\textit{Key words an phrases}. Lyapunov functions; Semiflows; Separable metric space.}The aim of this short note is to show how to construct a complete Lyapunov function of a semiflow by using a complete Lyapunov function of its time-one map. As a byproduct we assure the existence of complete Lyapunov functions for semiflows on separable metric spaces.
\end{abstract}


\section{Complete Lyapunov functions}

Let $\phi^t : X \to X$ be a semiflow on a metric space $X$ and denote by $\mathcal{R}_C(\phi^t)$ its chain recurrent set. A complete Lyapunov function for the semiflow $\phi^t$ is a continuous function $L : X \to [0,1]$ satisfying the following requirements:
\begin{enumerate}
\item The function $t \mapsto L(\phi^t(x))$ is constant for each $x \in \mathcal{R}_C(\phi^t)$ and is decreasing for each $x \in X - \mathcal{R}_C(\phi^t)$.

\item The set $L(\mathcal{R}_C(\phi^t))$ is nowhere dense in $[0,1]$.

\item If $c \in L(\mathcal{R}_C(\phi^t))$, then $L^{-1}(c)$ is a chain transitive component of $\mathcal{R}_C(\phi^t)$.
\end{enumerate}
A complete Lyapunov function for a continuous map $T : X \to X$ is a continuous function $L : X \to [0,1]$ satisfying the above requirements, changing the function $t \mapsto L(\phi^t(x))$ by the function $n \mapsto L(T^n(x))$ and changing the chain recurrent set of $\phi^t$ by the chain recurrent set of $T$, denoted by $\mathcal{R}_C(T)$.

The existence of a complete Lyapunov function for flows on compact metric spaces was originally proved by Conley in \cite{conley}. After this, in the reference \cite{franks}, Franks proved the existence of a complete Lyapunov function for maps on compact metric spaces, which was extended to locally compact and then to separable metric spaces by Hurley respectively in \cite{hurley1} and \cite{hurley3}. In some places (see e.g. \cite{alongi-nelson, norton, patrao}), these existence results are referred as the Fundamental Theorem of Dynamical Systems.

The aim of this short note is to show how to construct a complete Lyapunov function of a semiflow by using a complete Lyapunov function of its time-one map. As a byproduct we assure the existence of complete Lyapunov functions for semiflows on separable metric spaces.

\begin{teorema}
Let $\phi^t$ be a semiflow on $X$ and $\ell : X \to [0,1]$ be a complete Lyapunov function for its time-one map $\phi^1$. Thus the map $L : X \to [0,1]$ given by
\[
L(x) = \int_0^1 \ell(\phi^t(x)) \, dt
\]
is a complete Lyapunov function for the semiflow $\phi^t$. In particular, if $X$ is a separable metric space, a complete Lyapunov function always exists for any given semiflow.
\end{teorema}
\begin{profe}
We first note that $L(x)$ is well defined, for each $x \in X$, since the function $t \mapsto \ell(\phi^t(x))$ is continuous. It is straightforward to show that $L$ is continuous at any given $x \in X$. Indeed, given $\varepsilon > 0$, using the continuity of the semiflow and of $\ell$, for each $t \in [0,1]$, there exist open neighborhoods $U_t$ and $V_t$, respectively, of $t$ and of $x$ such that
\[
|\ell(\phi^s(y)) - \ell(\phi^t(x))| < \frac{\varepsilon}{2},
\]
for every $s \in U_t$ and $y \in V_t$. Since $\{U_t : t \in [0,1]\}$ is an open cover of the compact set $[0,1]$, there exists a finite subcover $\{U_{t_1}, \ldots, U_{t_n}\}$ and thus $V = V_{t_1} \cap \cdots \cap V_{t_n}$ is an open neighborhood of $x$. Hence, for every $t \in [0,1]$ and $y \in V$, we have that $t \in U_{t_k}$ and $y \in V_{t_k}$, for some $k$. Therefore, for every $t \in [0,1]$ and every $y \in V$, we have that
\[
|\ell(\phi^t(y)) - \ell(\phi^t(x))| \leq |\ell(\phi^t(y)) - \ell(\phi^{t_k}(x))| + |\ell(\phi^{t_k}(x)) - \ell(\phi^t(x))| < \varepsilon
\]
and thus
\[
|L(y) - L(x)| \leq \int_0^1 |\ell(\phi^t(y)) - \ell(\phi^t(x))| \, dt < \varepsilon,
\]
showing that $L$ is continuous at $x$.

Now we need to prove the three above requirements. We start by using reference \cite{hurley2} to remember that $\mathcal{R}_C(\phi^t) = \mathcal{R}_C(\phi^1)$, i.e., the chain recurrent sets of both the semiflow and its time-one map coincide. In order to verify the first requirement, for any $s \in (0,1]$, let us consider
\[
L(\phi^s(x)) = \int_0^1 \ell(\phi^t(\phi^s(x))) \, dt = \int_0^1 \ell(\phi^{t+s}(x)) \, dt =  \int_s^{1+s} \ell(\phi^{t}(x)) \, dt.
\]
Thus we have that
\begin{eqnarray*}
L(\phi^s(x)) & = & \int_s^{1} \ell(\phi^{t}(x)) \, dt + \int_1^{1+s} \ell(\phi^{t}(x)) \, dt \\
& = & \int_s^{1} \ell(\phi^{t}(x)) \, dt + \int_0^{s} \ell(\phi^{1+t}(x)) \, dt \\
& = & \int_s^{1} \ell(\phi^{t}(x)) \, dt + \int_0^{s} \ell(\phi^{1}(\phi^{t}(x)) \, dt
\end{eqnarray*}
If $x \in \mathcal{R}_C(\phi^t)$, then $\ell(\phi^{t}(x)) \in \mathcal{R}_C(\phi^1)$, implying that $\ell(\phi^{1}(\phi^{t}(x)) = \ell(\phi^{t}(x))$ and hence that
\[
L(\phi^s(x)) = \int_0^{1} \ell(\phi^{t}(x)) = L(x).
\]
On the other hand, if $x \in X - \mathcal{R}_C(\phi^t)$, then $\ell(\phi^{t}(x)) \in X - \mathcal{R}_C(\phi^1)$, implying that $\ell(\phi^{1}(\phi^{t}(x)) < \ell(\phi^{t}(x))$ and hence that
\[
L(\phi^s(x)) < \int_0^{1} \ell(\phi^{t}(x)) = L(x).
\]
The first requirement is thus immediate from the above results.

For the remaining requirements, we first note that the chain transitive components of $\mathcal{R}_C(\phi^1)$ and of $\mathcal{R}_C(\phi^t)$ coincide. Indeed each chain transitive component of $\mathcal{R}_C(\phi^1)$ is contained in a chain transitive component of $\mathcal{R}_C(\phi^t)$, since the transitivity by $\phi^1$ implies the transitivity by $\phi^t$. On the other hand, each chain transitive component of $\mathcal{R}_C(\phi^t)$ is contained in a chain transitive component of $\mathcal{R}_C(\phi^1)$, since the chain transitive components of $\mathcal{R}_C(\phi^1)$ are unions of connected components of $\mathcal{R}_C(\phi^1)$, while the chain transitive components of $\mathcal{R}_C(\phi^t)$ are exactly the connected components of $\mathcal{R}_C(\phi^t)$.

Now consider $c \in \ell(\mathcal{R}_C(\phi^1))$. Since $\ell^{-1}(c)$ is a chain transitive component of $\mathcal{R}_C(\phi^1)$, it is a chain transitive component of $\mathcal{R}_C(\phi^1)$ and thus is invariant by $\phi^t$. Hence, if $x \in \ell^{-1}(c)$, we have that $\phi^t(x) \in \ell^{-1}(c)$, which implies that
\[
L(x) = \int_0^1 \ell(\phi^t(x)) \, dt = c,
\]
showing that $\ell^{-1}(c) \subset L^{-1}(c)$. Since $\{ \ell^{-1}(c) : c \in \ell(\mathcal{R}_C(\phi^1)) \}$ is a partition of $\mathcal{R}_C(\phi^t)$, we have that $\ell(x) \in [0,1] - \ell(\mathcal{R}_C(\phi^1))$, for every $x \in X - \mathcal{R}_C(\phi^t)$. Using that $X - \mathcal{R}_C(\phi^t)$ is invariant by $\phi^t$ and using the continuity of the map $t \mapsto \ell(\phi^t(x))$, we get that the set $\{\ell(\phi^t(x)) : t \in [0,1]\}$ is a closed interval contained in $[0,1] - \ell(\mathcal{R}_C(\phi^1))$. This implies that $L(x) \in [0,1] - \ell(\mathcal{R}_C(\phi^1))$, for every $x \in X - \mathcal{R}_C(\phi^t)$, showing that $L^{-1}(c)$ is contained in $\mathcal{R}_C(\phi^t)$, for every $c \in \ell(\mathcal{R}_C(\phi^1))$. Since $\ell^{-1}(c) \subset L^{-1}(c)$ and since $\{ \ell^{-1}(c) : c \in \ell(\mathcal{R}_C(\phi^1)) \}$ is a partition of $\mathcal{R}_C(\phi^t)$, we get that $\ell^{-1}(c) = L^{-1}(c)$, for every $c \in \ell(\mathcal{R}_C(\phi^1))$. This implies that $L(\mathcal{R}_C(\phi^t)) = \ell(\mathcal{R}_C(\phi^1))$, showing that it is nowhere dense in $[0,1]$ and that $L^{-1}(c)$ is a chain transitive component of $\mathcal{R}_C(\phi^t)$, for every $c \in L(\mathcal{R}_C(\phi^t))$.

The last assertion follows by the first part, applying to the time-one map of the semiflow the main result of \cite{hurley3}, which assures the existence of a complete Lyapunov function for maps on separable metric spaces.

\end{profe}

\end{document}